\documentclass[10pt]{article}

\usepackage{rotating}
\usepackage{graphicx}
\usepackage{color}
\usepackage{placeins}
\usepackage{amssymb, amsmath, amsthm}

\usepackage{fancyhdr}
\fancypagestyle{preprint}{\fancyhf{}\fancyfoot[C]{Preprint submitted to AIP Conference Proceedings.}}

\begin{document}

\newcommand{\dx}{\Delta x}
\newcommand{\dt}{\Delta t}
\newcommand{\Oh}{{\cal O}}
\newcommand{\R}{\mathbb{R}}
\newcommand\tch[4]{\multicolumn{#1}{#2}%
   {\begin{tabular}[#3]{@{}#2@{}}%
    \ignorespaces#4\unskip
    \end{tabular}}}

\begin{center}{\Large Patankar-Type Runge-Kutta Schemes for Linear PDEs}\\[.2cm]

S.~Ortleb\footnote{Fachbereich Mathematik und Naturwissenschaften, Universit\"at Kassel, Heinrich Plett Str. 40, 34132 Kassel, Germany. E-mail: ortleb@mathematik.uni-kassel.de},
W. Hundsdorfer\footnote{CWI, Science Park 123, Amsterdam, The Netherlands. E-mail: willem.hundsdorfer@cwi.nl}

\end{center}

\begin{abstract}
We study the local discretization error of Patankar-type Runge-Kutta methods applied to semi-discrete PDEs. For a known two-stage Patankar-type scheme the local error in PDE sense for linear advection or diffusion is shown to be of the maximal order ${\cal O}(\Delta t^3)$ for sufficiently smooth and positive exact solutions. However, in a test case mimicking a wetting-drying situation as in the context of shallow water flows, this scheme yields large errors in the drying region. A more realistic approximation is obtained by a modification of the Patankar approach incorporating an explicit testing stage into the implicit trapezoidal rule.
\end{abstract}

\thispagestyle{preprint}

\section{Quasi-linear equations in production-destruction form arising from PDEs}

In the context of geobiochemical models, so-called production-destruction equations are frequently encountered. These models describe the time-evolution of non-negative quantities and often take into account some type of mass conservation. The underlying ODE systems describing the time-evolution of non-negative quantities $u(t)$ can usually be written in the form
$u'_i = \Big( \sum_{j=1}^m p_{ij}(u)\cdot u_j \Big) - q_i(u)\cdot u_i$
for $i=1,2,\ldots,m$,
with production terms $p_{ij}(u) u_j$ and destruction terms $q_i(u) u_i$
such that $p_{ij}(v), \,q_i(v) \,\ge\, 0$ for all $v\in\R^m$ and 
$i,j=1,2,\ldots,m$. Usually, we also have $p_{ii}(v) = 0$ for all $v\in\R^m$ and some mass conservation property such as $\sum_{j=1}^mp_{ji}(u) = q_i(u)$.
In vector form, we write
\begin{eqnarray}
\label{eq:PQ1}
u' \,=\, P(u)\,u \,-\, Q(u)\,u,
\end{eqnarray}
with matrix-valued functions $P, Q : \R^m \rightarrow \R^{m\times m}$
such that
$P(v) = \big(p_{ij}(v)\big) \ge 0, \
Q(v) = \mbox{diag}\big(q_i(v)\big)\ge0 $
for all $v\in\R^m$.
Setting $A(u) = P(u) - Q(u)$, this is written more shortly in the standard quasi-linear form $u' \,=\, A(u)\,u$.

While many interesting biochemical reactions fit into this framework, it also includes certain space-discretized partial differential equations, e.g. the heat equation discretized by second-order differences and the first-order upwind-discretized advection equation. In the context of shallow water flows discretized by the discontinuous Galerkin method, a production-destruction approach as in \cite{meister_ortleb:14} guarantees non-negativity of the water height for any time step size while still preserving conservativity. In that work, the production-destruction equations where specifically formulated in order to account for the production and destruction terms which influence the cell-wise water volume.

Generally, numerical methods discretizing (\ref{eq:PQ1}) are supposed to be positivity preserving, conservative and of sufficiently high order. While positivity preservation and conservativity may be directly carried over from the context of ODEs to that of PDEs, the issue of consistency and convergence is more subtle for PDEs. In this work, we hence take a closer look at the local discretization error of Patankar-type methods applied to systems arising from linear PDEs.

\section{The Patankar-Euler method and its modification}

The forward Euler method applied to (\ref{eq:PQ1}) will obviously be
positivity preserving if  we have $I - \dt\,Q(u^n) \ge 0$, but this requires
a very severe time step restriction on $\dt$ for stiff systems.
To avoid this, a variant was proposed by Patankar \cite{Patankar:1980}, originally in the context of source terms in heat transfer. This method given by
\begin{eqnarray}
\label{eq:EulPata}
u^{n+1} \,=\, u^n + \dt\,P(u^n) u^n - \dt\,Q(u^n) u^{n+1} \,
\end{eqnarray}
is unconditionally positivity preserving but not mass conserving. 
In addition, while (\ref{eq:EulPata}) is of order one in the ODE 
sense, consistency is lost for stiff problems such as the discretized heat equation. In fact, the semi-discrete 1D heat equation
$u'_i(t) = \frac{1}{\dx^2} \big(u_{i-1}(t) - 2 u_i(t) + u_{i+1}(t) \big)$ for
$i=1,2,\ldots,m$,
with spatial periodicity, i.e. $u_0(t) = u_m(t)$ and $u_{m+1}(t) = u_1(t)$ fits in the form (\ref{eq:PQ1}) with diagonal destruction matrix
$Q(u) = 2 \dx^{-2} I$. 
The Patankar-Euler scheme (\ref{eq:EulPata}), written out per component, now reads
$u_i^{n+1} \,=\, u_i^n + \frac{\dt}{\dx^2} \big(u_{i-1}^n - 2 u_i^{n+1} 
+ u_{i+1}^n \big)$.
This scheme is unconditionally positivity preserving as well as unconditionally contractive in the maximum norm. However, inserting exact PDE solution values
in the scheme, we obtain $u(x_i,t^{n+1}) \,=\, u(x_i,t^n) + \frac{\dt}{\dx^2} \big( u(x_{i-1},t^n) - 2 u(x_i,t^{n+1})
+ u(x_{i+1},t^n) \big) + \dt \rho_i^n$. Hence, Taylor development shows that
for small $\dt$ and $\dx$ the leading term in these local truncation errors
is given by
$\rho_i^n \,=\, \frac{2 \dt}{\dx^2} u_t(x_i,t^n) + \Oh\Big(\frac{\dt^2}{\dx^2}\Big)$.
It follows that the scheme will only be convergent in case of a very severe time step restriction of $\dt/\dx^2 \rightarrow 0$.

In order to obtain an unconditionally positive and additionally mass conservative scheme for production-destruction equations, the modification
$u^{n+1} \,=\, u^n + \dt\,P(u^n) u^{n+1} - \dt\,Q(u^n) u^{n+1} \,=\, u^n + \dt\,A(u^n) u^{n+1}$ has been proposed in \cite{Burchard:2003}.
For linear problems $u'= A u$ with constant matrix $A$, such as the linear heat equation, this modified method now reduces to the implicit Euler method, so consistency in PDE sense is not a problem there.
In \cite{Burchard:2003}, also a second-order
method has been proposed. We will refer to this method as mPaRK2. This scheme does not fit directly in
the vector production-loss formulation and thus has to be written  
per component, starting with the quasi-linear form
$u'_i \,=\, \sum_{j=1}^m a_{ij}(u) u_j, i=1,2,\ldots,m$.
The mPaRK2 method is then based on the trapezoidal rule with an Euler-type prediction to provide the
internal stage value $v_i^{n+1} \approx u_i(t^{n+1})$ and reads
\begin{eqnarray}
\label{eq:MPRK}
v_i^{n+1} \, = \, \displaystyle
u_i^n \,+\, \dt \sum_{j} a_{ij}(u^n) v_j^{n+1} \,,
\quad
u_i^{n+1} &=& 
u_i^n \,+\, \frac{1}{2} \dt \sum_{j} \Big(
a_{ij}(u^n) \frac{u_j^n}{v_j^{n+1}}\,u_j^{n+1}
+ a_{ij}(v^{n+1}) u_j^{n+1} \Big) \,.
\end{eqnarray}
As shown in \cite{Burchard:2003}, this scheme is unconditionally 
positivity preserving and mass conserving, and the order is two in the 
ODE sense. 
However, it is unknown whether there will be order reduction for stiff problems,
in particular for semi-discrete problems obtained from PDEs after space
discretization. Regarding the local discretization error, consistency of ${\cal O}(\Delta t^3)$ can be proven for sufficiently smooth exact solutions. This is dealt with in the next section.

\section{Error recursions and numerical results}

We will study error recursions for the mPaRK2 method applied to linear problems with constant coefficients.
These are naturally non-linear for this method, even for linear equations.
For a linear problem $u'(t) = A u(t)$ with $A = (a_{ij}) \in\R^{m\times m}$, 
we will first write (\ref{eq:MPRK}) in vector form by introducing the 
diagonal matrix $W^n = \mbox{diag}(u_i^n/v_i^{n+1})$. Then (\ref{eq:MPRK})
can be written compactly as
\begin{eqnarray}
\label{eq:MPRKv}
v^{n+1} \, = \, \displaystyle
u^n \,+\, \dt \,A v^{n+1} \,,
\qquad
u^{n+1} \, = \,
u^n \,+\, \frac{1}{2} \dt\,A (W^n + I) u^{n+1} \,. 
\end{eqnarray}
Along with this, we also consider the scheme with the exact solution inserted,
\begin{eqnarray}
\label{eq:pMPRKv}
\bar{v}^{n+1} \, = \, \displaystyle
u(t^n) \,+\, \dt \,A \bar{v}^{n+1} \,,
\qquad
u(t^{n+1}) \, = \,
u(t^n) \,+\, \frac{1}{2} \dt\,A (\bar{W}^n + I) u(t^{n+1}) 
\,+\, \rho^n \,,
\end{eqnarray}
where $\bar{W}^n = \mbox{diag}\big(u_i(t^n)/\bar{v}_i^{n+1}\big)$ and $\rho^n = (\rho_i^n) \in \R^m$. 
Subtraction of (\ref{eq:MPRKv}) from (\ref{eq:pMPRKv}) gives
a recursion for the global discretization errors $e^n = u(t^n) - u^n$ of the form $e^{n+1} \,=\, R^ n e^n + d^n$,
with amplification matrix and local errors given by
\begin{eqnarray*}
R^n \,=\, 
\Big( I - \frac{1}{2}\dt A\big(\bar{W}^n + I\big) \Big)^{-1} 
\Big( I + \frac{1}{2}\dt A G^n \Big)\,, \quad
d^n \,=\, \Big( I - \frac{1}{2}\dt A\big(\bar{W}^n + I\big) \Big)^{-1} \rho^n
 \,,
\end{eqnarray*}
with the matrix $G^n \in \R^{m\times m}$ given by 
$G^n \,=\, \mbox{diag}({{u}_i^{n+1}}/{\bar{v}_i^{n+1}} )
- \mbox{diag}((u_i^n {u}_i^{n+1})/
(\bar{v}_i^{n+1} v_i^{n+1}) )
(I-\dt A)^{-1}$.
The difference between $\rho^n$ and its counterpart resulting from the implicit trapezoidal rule can be determined from
\begin{eqnarray}
\begin{array}{ccl}
  \rho^n = u(t^{n+1}) - u(t^n) - \frac{1}{2}\dt A\left(u(t^{n}) + u(t^{n+1})\right) + \frac{1}{2}\dt A\left(u(t^{n}) - \bar{W}^nu(t^{n+1})\right). 
\end{array}
\end{eqnarray}
Thus, the term
$\tilde{\rho}^n \, = \,\frac{1}{2}\dt A\left(u(t^{n}) - \bar{W}^{n}u(t^{n+1})\right) \, = \, \frac{1}{2}\dt A\,\mbox{diag}\left(\bar{v}_i^{n+1} - u_i(t^{n+1})\right)\left(\mbox{diag}\left(\bar{v}_i^{n+1}\right)\right)^{-1}u(t^n) = \frac{1}{2}\dt A\, D_1\, D_2\, u(t^n)$
represents the difference in local errors between the implicit trapezoidal rule and the mPaRK2 scheme. For the diagonal matrices, we have
$  D_1 \, = \, \mbox{diag}\left(\left((I-\Delta t A)^{-1}-e^{\Delta t A}\right)u(t^{n})\right) = {\cal O}(\Delta t^2)$ and
$  D_2 \, = \, \left(\mbox{diag}\left((I-\Delta t A)^{-1}u(t^{n})\right)\right)^{-1}$.
In addition, if we assume $u(t^n)>0$ then $D_2$ is bounded for $\Delta t \rightarrow 0$, i.e.  $D_2 = {\cal O}(1)$. 

\begin{paragraph}{Semi-discrete linear advection and linear diffusion}
A reasonable assumption in the case of the semi-discrete linear advection with $A = \frac{1}{\dx}\,\mbox{tridiag}[1\, -1 \;\; 0]$ is a time-step choice such that $\Delta t A  = {\cal O}(1)$. Then we have
\begin{eqnarray}
\begin{array}{ccl}
  \tilde{\rho}^n & = &\frac{1}{2}\underbrace{\dt A}_{{\cal O}(1)} \left(\underbrace{\left((I-\Delta t A)^{-1}-e^{\Delta t A}\right)}_{{\cal O}(\Delta t^2)} + \underbrace{\mbox{diag}\left(\frac{u_i(t^{n})-\bar{v}_i^{n+1}}{\bar{v}_i^{n+1}}\right)}_{{\cal O}(\Delta t)} \underbrace{\left((I-\Delta t A)^{-1}-e^{\Delta t A}\right)}_{{\cal O}(\Delta t^2)}\right) u(t^n)\\
  & = & \frac{1}{2}\dt A\left((I-\Delta t A)^{-1}-e^{\Delta t A}\right)u(t^n) + {\cal O}(\Delta t^3).
\end{array}
\end{eqnarray}
Due to the smoothness of the implicit Euler scheme, it holds that $A\left((I-\Delta t A)^{-1}-e^{\Delta t A}\right)u(t^n) = {\cal O}(\Delta t^2)$. Hence, the local discretization error of the mPaRK2 method applied to the semi-discrete linear advection equation with positive initial data is bounded by $\rho^n = {\cal O}(\Delta t^3)$.

Additional assumptions on the smoothness of the initial data are necessary for the semi-discrete linear diffusion equation with $A = \frac{1}{\dx^2}\,\mbox{tridiag}[1\, -2 \;\; 1]$ as shown in Fig. \ref{fig:smooth} and Table \ref{tab:loc_err}.
As the diagonal matrix $D_1$ in the definition of $\tilde{\rho}^n$ 
can be bounded by $D_1  =  \mbox{diag}\left(\frac{1}{2}(\Delta t A)^{2}u(t^n)\right) + {\cal O}((\Delta t A)^3)$,
a smoothness condition on the exact solution of the type
\begin{eqnarray}\label{eq:cond_smooth}
S:=  A\,\mbox{diag}\left((A^2u(t^{n}))_i\right) \left(\mbox{diag}((I-\Delta t A)^{-1}u(t^{n}))\right)^{-1}u(t^{n}) = {\cal O}(1)
\end{eqnarray}
guarantees $\tilde{\rho}^n = {\cal O}(\Delta t^3)$ and hence a local error of third order.
The left part of Fig. \ref{fig:smooth} depicts the situation for an initial solution $u_0 = 0.1 + \sin^2(2\pi x)$ which satisfy the smoothness condition. Here, the quantities $u_i/v_i$ approximate 1 for $\Delta t\rightarrow 0$ with $\Delta t = {\cal O}(\Delta x)$. Consequently, the difference between the quantities $A^2u_0\approx u_{0,xxxx}$ and $\mbox{diag}\left(u_i/v_i\right)A^2u_0$ is small leading to basically constant smoothness indicators $\|S\|_2$ as shown in Table \ref{tab:loc_err}.
This Table also lists the local error in the first mPaRK2 step. In accordance with the designed order of convergence, this local error behaves as ${\cal O}(\Delta t^3)$. On the other hand, for an initial solution of $u_0 = \sin^2(2\pi x)$, Table \ref{tab:loc_err} shows an order reduction to about ${\cal O}(\Delta t^{2.3})$ in addition to an increasing value of the indicator $\|S\|_2$. As depicted in Fig. \ref{fig:smooth}, this behavior is due the fact that for $u_i = 0$, we also have $u_i/v_i = 0$. Values at nearby grid points will tend to 1 for $\Delta t\rightarrow 0$ while zeros of $u_i/v_i$ remain unchanged. This leads to the boundary layer effect for $u_i/v_i = 0$ visible in the right part of Fig. \ref{fig:smooth} as well as a locally large difference in curvature between $A^2u_0$ and $\mbox{diag}\left(u_i/v_i\right)A^2u_0$.
\begin{table}[h]
\caption{Effect of the smoothness condition (\ref{eq:cond_smooth}) on the error of consistency for the mPaRK2 scheme.}
\label{tab:loc_err}
\tabcolsep7pt\begin{tabular}{l|ccc|ccc}
\hline
$m$  & \tch{3}{c}{b}{$u_0 = 0.1 + \sin^2(2\pi x)$} & \tch{3}{c}{b}{ $u_0 = \sin^2(2\pi x)$} \\
  & $L^2$ loc. err. & $L^2$ consistency  & $\|S\|_2$  & $L^2$ loc. err.  & $L^2$ consistency  &  $\|S\|_2$   \\
\hline
$40$  & 0.00177  &  & 1.70e+06 & 0.00218 &     & 2.55e+06 \\
$80$  & 0.00036  & 2.31 & 1.65e+06 & 0.00054 &   2.02  & 3.40e+06 \\
$160$ & 5.74e-05 & 2.64 & 1.53e+06 & 0.00012 &  2.20  & 4.89e+06 \\
$320$ & 8.13e-06 & 2.82 & 1.44e+06 & 2.45e-05 & 2.25  & 7.52e+06 \\
$640$ & 1.08e-06 & 2.91 & 1.39e+06 & 5.12e-06 & 2.26  & 1.21e+07 \\
$1280$& 1.40e-07 & 2.95 & 1.38e+06 & 1.07e-06 & 2.26 & 1.98e+07 \\
$2560$& 1.78e-08 & 2.97 & 1.81e+06 & 2.24e-07 & 2.25 & 3.30e+07 \\
\hline
\end{tabular}
\end{table}
\begin{figure}[h]
\begin{minipage}{7.5cm}
\includegraphics[width=7.5cm,height=4.5cm]{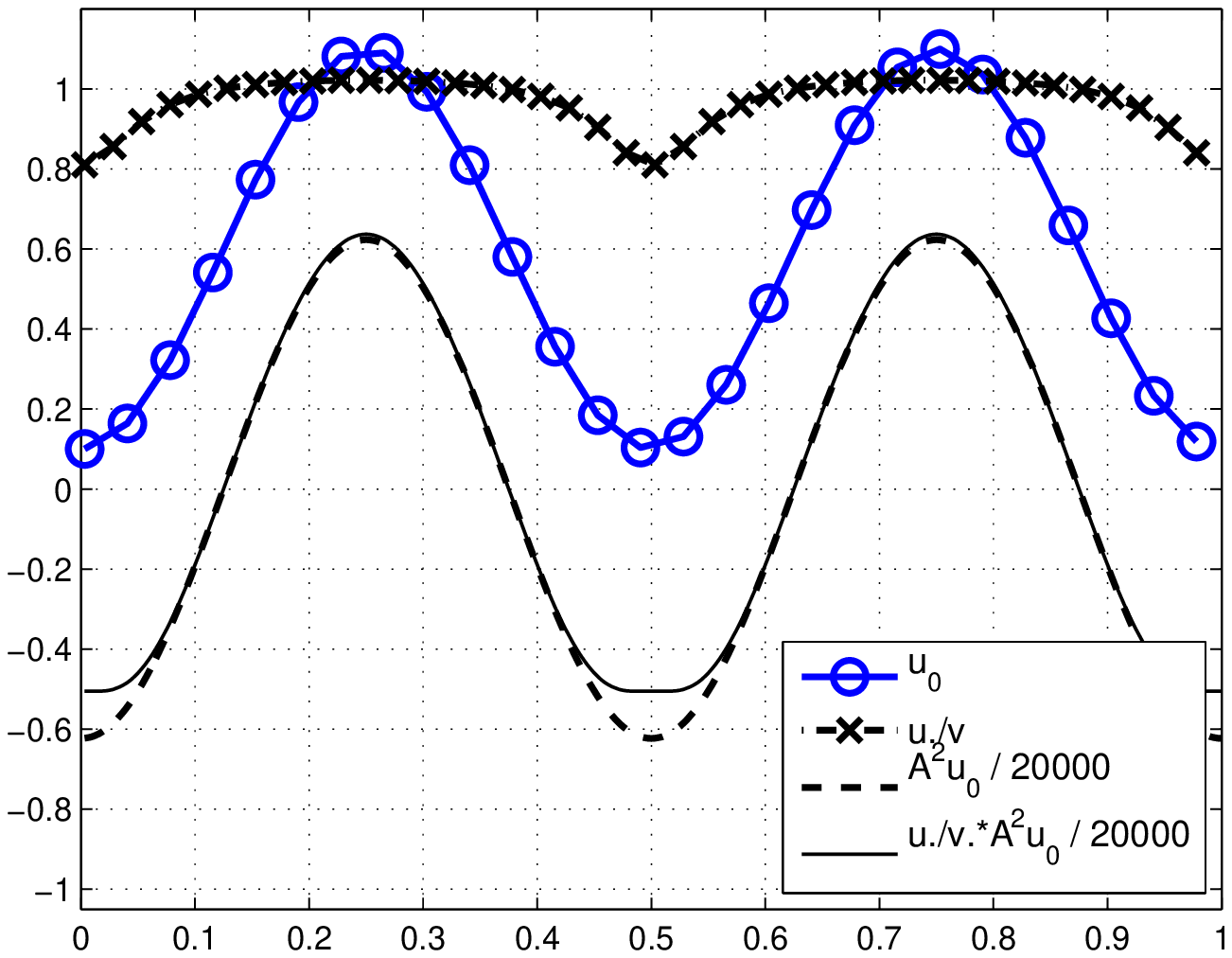}
\end{minipage}
\hspace*{.2cm}
\begin{minipage}{7.5cm}
\includegraphics[width=7.5cm,height=4.5cm]{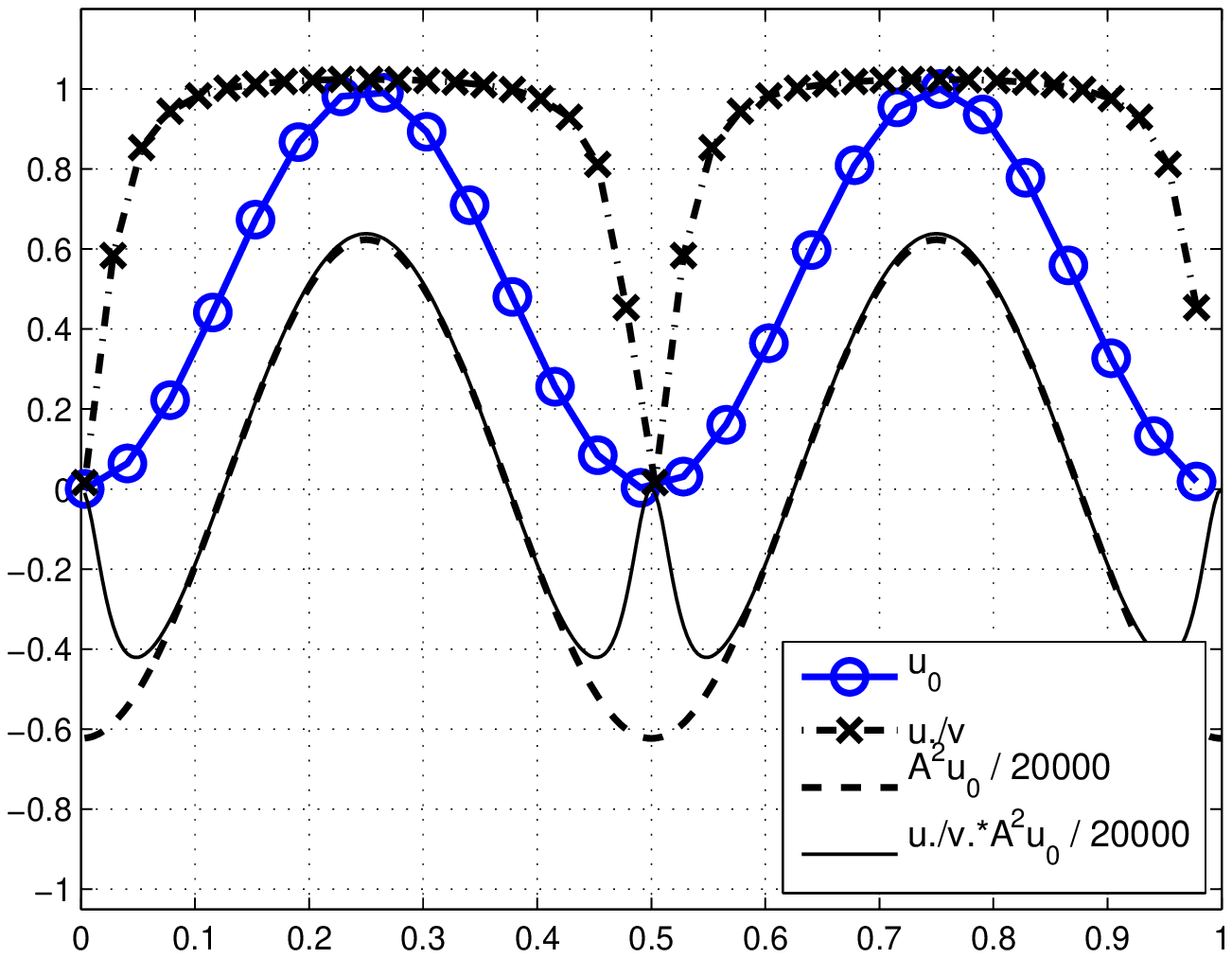}
\end{minipage}
\caption{Initial conditions for the semi-discrete linear heat equation which satisfy (left) or do not satisfy (right) condition (\ref{eq:cond_smooth}).}
 \label{fig:smooth}
\end{figure}
\end{paragraph}

So far, these investigations show a local discretization error of order ${\cal O}(\Delta t^3)$ for sufficiently smooth and positive solutions. This positivity requirement includes thin-layer approaches for the shallow water equations, where a thin film of water is retained also in regions marked as dry. However, we should remark that for a full convergence analysis, stability has to be proven as well. This necessitates boundedness of products of amplification matrices $R^n$ which seems to be quite difficult to prove due to the non-linearity of the method.

\begin{figure}[h]
\begin{minipage}{7.5cm}
\includegraphics[width=7.5cm,height=4.6cm]{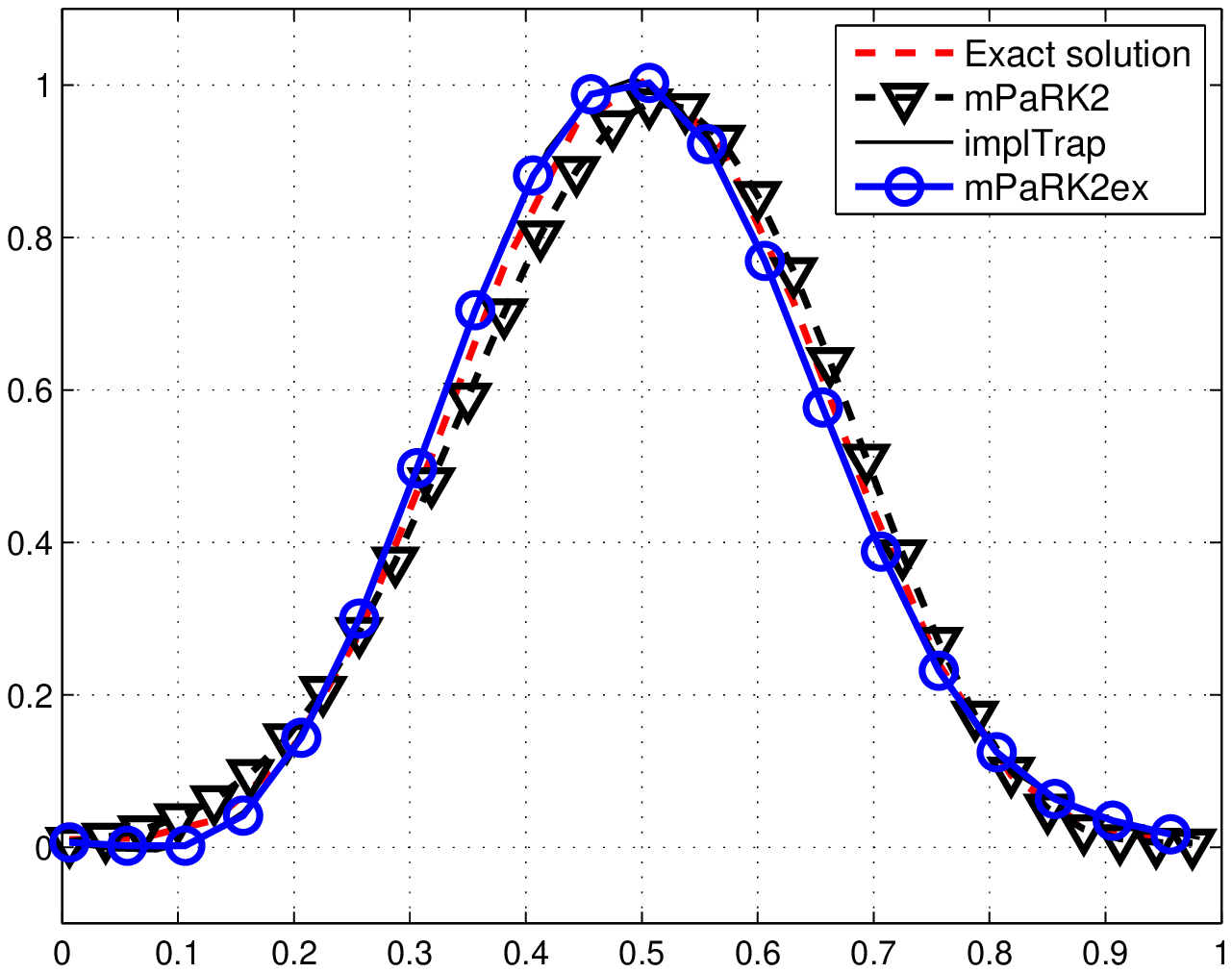}
\end{minipage}
\hspace*{.2cm}
\begin{minipage}{7.5cm}
\includegraphics[width=7.5cm,height=4.6cm]{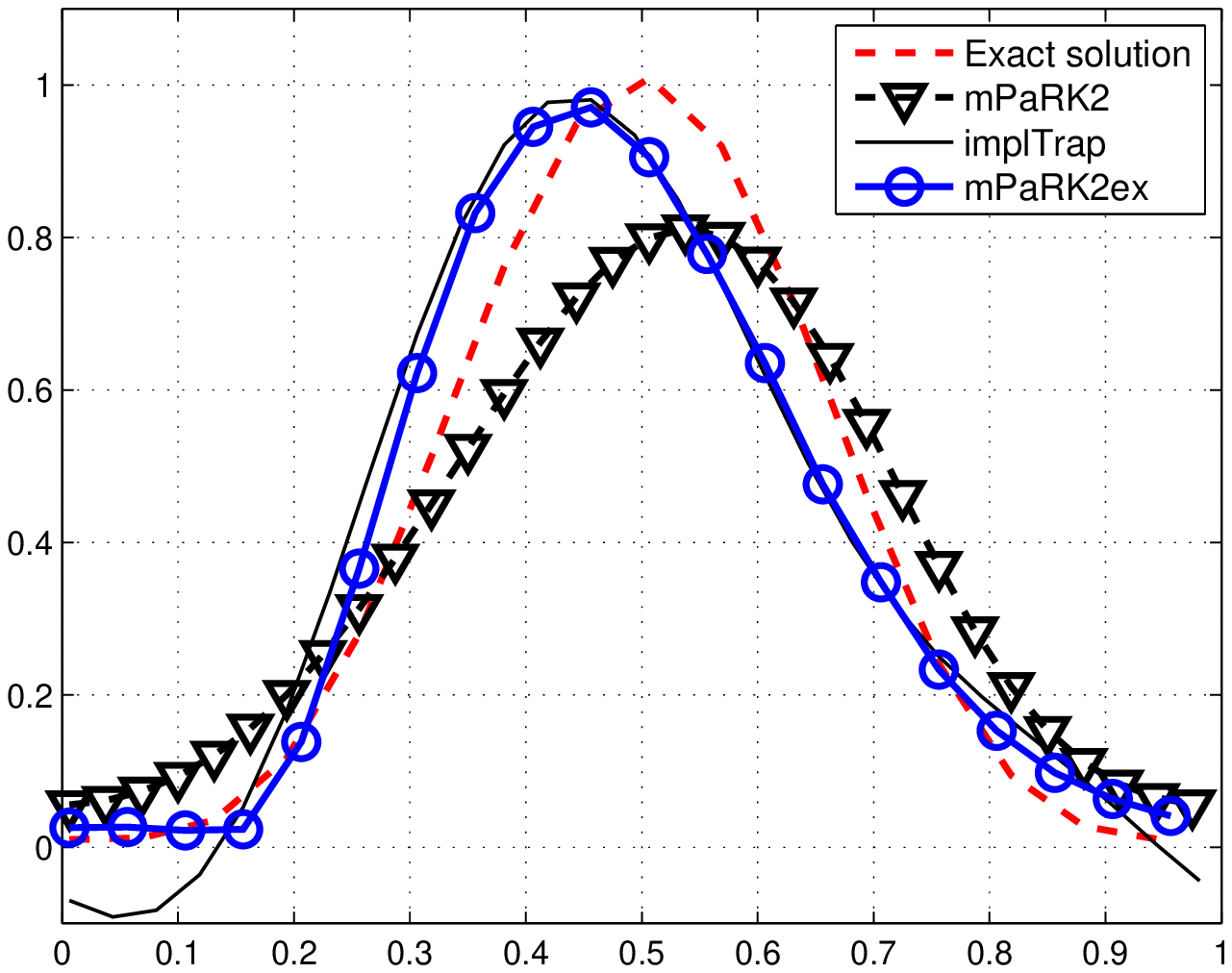}
\end{minipage}
\caption{Linear transport of a wave using a Courant number of $4$ (left) and $10$ (right). Comparison of mPaRK2 to mPaRK2ex.}
 \label{fig:patSWE}
\end{figure}

Finally, we also consider a modification to the mPaRK2 scheme which follows more closely the approach in \cite{meister_ortleb:14}. This modification is based on a direct correction of the explicit part of the implicit trapezoidal rule and reads as
\begin{eqnarray}
\label{eq:mPaRK2ex}
v^{n+1/2} \, = \, \displaystyle
u^n \,+\, \frac{\dt}{2} \,A u^n \,,
\quad
u^{n+1/2} &=& 
u^n \,+\, \frac{1}{2} \dt A \tilde{W}\,u^{n+1/2}\,,
\quad
u^{n+1} \, = \,
u^{n+1/2} \,+\, \frac{1}{2} \dt\,A u^{n+1} \,,
\end{eqnarray}
with $\tilde{W} = \mbox{diag}\left(u_i^{n} / \tilde{v}_i^{n+1/2}\right)$ determined by the correction $\tilde{v}^{n+1/2}$ to the quantity $v^{n+1/2}$ which may have negative components. More precisely, $\tilde{v}^{n+1/2}$ is given by $\tilde{v}_i^{n+1/2} = v_i^{n+1/2}$ if $v_i^{n+1/2}>0$ and $\tilde{v}_i^{n+1/2} = u_i^{n}$ otherwise. We will denote this scheme by mPaRK2ex. Due to this switch in case of vanishing components, we cannot expect an overall second order of convergence as the update reduces to two steps of the implicit Euler scheme if $v^{n+1/2}=0$. However, for a test case of an advected wave mimicking wetting and drying, i.e. advection of the initial condition $u_0 = 0.01 + \sin^4(\pi x)$, this method behaves much better than mPaRK2 as shown in Fig. \ref{fig:patSWE}.

A comparison of the Patankar-type schemes is carried out for the upwind-discretized linear advection on $160$ grid points using spatial periodicity up to a final time of $T = 2$. As shown on the left of Fig. \ref{fig:patSWE}, using a time step of $\Delta t = 0.025$ corresponding to a Courant number of $4$ does not exhibit significant differences of the schemes mPaRK2 and mPaRK2ex, also in comparison to the implicit trapezoidal rule. However, a larger time step of $\Delta t = 0.0625$ corresponding to a Courant number of $10$ shows the drawback of mPaRK2 on the right part of Fig. \ref{fig:patSWE}. While mPaRK2 does not account for the vanishing solution in the interval $[0, 0.15]$ and the implicit trapezoidal rule clearly yields negative values, the modified scheme mPaRK2ex seems to combine the best features of both methods. The solution is non-negative in the whole computational domain and very accurate in the almost dry regions. Hence, this method seems quite promising and should be further investigated, in particular with respect to its stability.

\section*{Acknowledgements}
S. Ortleb gratefully acknowledges the hospitality of the Centrum Wiskunde \& Informatica where this research was partly carried out during her visit in the period of February/March 2016.



\begin{thebibliography}{}
\bibitem{Burchard:2003} Burchard, H., Deleersnijder, E., Meister, A.: A High-order Conservative Patankar-type Discretisation for Stiff Systems of Production-destruction Equations. Appl. Numer. Math. \textbf{47}, 1--30 (2003)
\bibitem{Patankar:1980} Patankar, S. V.: \emph{Numerical heat transfer and fluid flow}. Series in computational methods in mechanics and thermal sciences, Hemisphere Pub. Corp. New York, Washington, 1980.
\bibitem{meister_ortleb:14} Meister, A., Ortleb, S.: On unconditionally positive implicit time integration for the DG scheme applied to shallow water flows, International Journal for Numerical Methods in Fluids \textbf{76}, 69--94 (2014)
  \end{thebibliography}
\end{document}